\documentclass[leqno,12pt]{article}
\usepackage{amsmath,amssymb,latexsym}
\usepackage{graphics}
\usepackage{amsmath,amssymb,latexsym,euscript}
\usepackage[T1]{fontenc}
\titlepage
\newcommand{\bea}{\begin{eqnarray*}}
\newcommand{\eea}{\end{eqnarray*} }
\newcommand{\be}{\begin{eqnarray}}
\newcommand{\ee}{\end{eqnarray} }
\newcommand{\beq}{\begin{equation}}
\newcommand{\eeq}{\end{equation} }
\newcommand{\CCp}{\mathcal{P}}
\newcommand{\CCs}{\mathcal{S}}

\newcommand{\BBone}{\mbox{1}\kern-.25em \mbox{I}}

\newtheorem{theorem}{Theorem}[section]
\newtheorem{lemma}{Lemma}[section]
\newtheorem{proposition}{Proposition}[section]
\newtheorem{corollary}{Corollary}[section]
\newtheorem{definition}{Definition}[section]
\newtheorem{example}{Example}[section]
\newtheorem{algorithm}{algorithm}[section]
\newtheorem{remark}{Remark}[section]
\newtheorem{statement}{Statement}[section]
\newcommand{\bt}{\begin{theorem}}
\newcommand{\et}{\end{theorem}}
\newcommand{\bl}{\begin{lemma}}
\newcommand{\el}{\end{lemma}}
\newcommand{\bp}{\begin{proposition}}
\newcommand{\ep}{\end{proposition}}
\newcommand{\bc}{\begin{corollary}}
\newcommand{\ec}{\end{corollary}}
\newcommand{\btb}{\begin{table}[hbp]}
\newcommand{\btbh}{\begin{table}[hhh]}
\newcommand{\etb}{\end{table}}
\newcommand{\bfg}{\begin{figure}[hbp]}
\newcommand{\bfgh}{\begin{figure}[hhh]}
\newcommand{\efg}{\end{figure}}
\newcommand{\bd}{\begin{definition}\rm}
\newcommand{\ed}{\end{definition}}
\newcommand{\bex}{\begin{example}\rm}
\newcommand{\eex}{\end{example}}
\newcommand{\ba}{\begin{algorithm}\rm}
\newcommand{\ea}{\end{algorithm}}
\newcommand{\br}{\begin{remark}\rm}
\newcommand{\er}{\end{remark}}
\newcommand{\bs}{\begin{statement}}
\newcommand{\es}{\end{statement}}

\def\er{\mathbb{R}}
\begin{document}

\title{A note on random orthogonal polynomials\\ on a compact interval}

\author{Melanie Birke \\
Ruhr-Universit\"at Bochum \\
Fakult\"at f\"ur Mathematik \\
44780 Bochum, Germany \\
{\small e-mail: melanie.birke@ruhr-uni-bochum.de }\\
\and
Holger Dette\\
Ruhr-Universit\"at Bochum \\
Fakult\"at f\"ur Mathematik \\
44780 Bochum, Germany \\
{\small email: holger.dette@ruhr-uni-bochum.de}\\
}

\maketitle

\begin{abstract} \noindent
We consider a uniform distribution on the set $\mathcal{M}_k$ of moments of order $k \in \mathbb{N}$ corresponding to probability measures on
the interval $[0,1]$. To each (random) vector of moments in $\mathcal{M}_{2n-1}$ we consider the corresponding  uniquely determined
 monic (random) orthogonal polynomial of degree $n$ and study the asymptotic properties of its roots if $n \to \infty$.
\end{abstract}

\medskip


\section{ Introduction } 
\def\theequation{1.\arabic{equation}}
\setcounter{equation}{0}

For a probability measure $\eta$ on the interval $[0,1]$ let \be \label{1.1} c_k = c_k (\eta) = \int^1_0 x^k \eta (dx) \qquad k = 0,1,2,\dots
\ee denote the corresponding moments and consider the moment space (of order $n$) defined by \be \label{1.2} \mathcal{M}_n = \{  (c_1 (\eta),
\dots,  c_n (\eta))^T \ | \quad \eta \  \mbox{is a probability measure on the interval}  \ [0,1]  \ \} \ . \ee The set $\mathcal{M}_n$  is a
very small subset of $\mathbb{R}^n$ with volume proportional to $2^{-n^2}$ and has been studied extensively in the literature [see e.g.\ Karlin
and Shapeley (1953), Skibinsky (1967, 1968, 1969), Chang, Kemperman and Studden (1993) or Gamboa and Lozada-Chang (2004)].  In order to
understand the structure of the moment space Chang, Kemperman and Studden (1993) assigned a uniform distribution over $\mathcal{M}_n$ and
studied the asymptotic properties of the first $k$ components of the random vector \be \label{1.3} (C_{1,n}, \dots, C_{n,n})^T \sim \mathcal{U}
(\mathcal{M}_n) \ee if $n \to \infty$. In particular, they showed that an appropriately standardized version of the vector $(C_{1,n} \dots,
C_{k,n})^T$ converges weakly to a centered normal distribution, that is \be \label{1.4} \sqrt{n} \ \{  (C_{1,n}, \dots, C_{k,n})^T - (c^0_1,
\dots, c^0_k)^T \} \rightarrow \mathcal{N} (0, \Sigma), \ee where $c^0_k$ denotes the $k$th moment of the arcsine distribution defined by \be
\label{1.5} c^0_k = \frac {1}{\pi} \int^1_0 \frac {x^k}{\sqrt{x(1-x)}} \ dx = \frac {1}{2^{2k}}  {2k \choose k}  \ , \qquad k = 0,1,\dots \ ,
\ee and the matrix $\Sigma \in \mathbb{R}^{k \times k}$ is given by \be \label{1.6} \Sigma = \frac {1}{2} \left(c^0_{i+j} - c^0_i
c^0_j\right)^k_{i,j=1} \ . \ee A large deviation principle for the random moment vector $(C_{1,n} \dots, C_{k,n})^T$ was derived in Gamboa and
Lozada-Chang (2004), while Dette and Gamboa (2007) investigated the asymptotic properties of a moment range process.

It is the purpose of the present paper to provide further insight in the probabilistic properties of quantities associated to random moment
sequences. In particular, by the definition $\mathcal{L}(x^k) = C_{k,2n-1} \quad (k = 1, \dots, 2n-1)$ each random vector $(C_{1,2n-1}, \dots,
C_{2n-1, 2n-1})^T \sim \mathcal{U} (\mathcal{M}_{2n-1})$ defines a (random) moment functional and a corresponding sequence of monic random
orthogonal polynomials $P_{0,n} (x), \dots, P_{n,n} (x)$ satisfying \be \label{1.7} \mathcal{L} (P_{k,n} \ P_{l,n}) = 0 \quad \mbox{if} \quad k
\neq l \ ; \ k,l \in \{ 0, \dots, n \} \ee [see Chihara (1978)]. These polynomials are carefully introduced in Section 2, where we also state
some non-standard results regarding moment theory.
 In Section 3 we derive the asymptotic properties of the (random)
roots $X_{1,n}, \dots, X_{n,n}$ of the polynomial $P_{n,n} (x)$ associated with the random moment sequence $(C_{1,2n-1}, \dots, C_{2n-1,2n-1})
\sim \mathcal{U}(\mathcal{M}_{2n-1})$. In particular, it is shown that the empirical distribution function of the (random) roots $ X_{1,n},
\dots, X_{n,n} $ converges almost surely to the arcsine distribution if $n \to \infty$ and that an appropriate standardization of the vector
$X_{1,n}, \dots, X_{k,n}$ is asymptotically normal distributed, where the roots of the Chebyshev polynomial of the first kind are used for the
centering, and the rate of convergence is $1/\sqrt{n}$.

\section{ Random coefficients in a three term recurrence relation } 
\def\theequation{2.\arabic{equation}}
\setcounter{equation}{0}

The set ${\cal M}_{n}$  defined by (\ref{1.2}) is a very small compact subset of the unit cube $[0,1]^{n}$ with nonempty interior and volume
\begin{eqnarray} \label{2.0}
\mbox{Vol}({\cal M}_{n}) = \prod^{2n}_{k=1} \frac{\Gamma (k) \Gamma(k)}{\Gamma(2k)} \approx c \cdot 2^{-n^2}
\end{eqnarray}
[see Karlin and Shapely (1953)]. The interior of $\mathcal{M}_n$ is  denoted by $\mathcal{M}^0_n$ throughout this paper. It is well known that
there exist an infinite number of probability measures on the interval $[0,1]$ with moments up to the order $2n-1$ given by $(c_1, \dots,
c_{2n-1})^T \in \mathcal{M}^0_{2n-1}$ [see Dette and Studden (1997)]. Moreover, there exists a unique measure $\eta^-$ supported on exactly $n$
points in the open  interval (0,1) such that $c_i(\eta^-) = \int^\prime_0 x^i d \eta^- (x) = c_i (i=1, \dots, 2n-1)$ and such that the point
$(c_1 (\eta^-), \dots, c_{2n} (\eta^-))^T \in
\partial  \mathcal{M}_{2n}$ where $\partial \mathcal{M}_{2n}$ denotes the boundary of the set $\mathcal{M}_{2n}$.
The measure $\eta^-$ is called lower principal representation of the point $(c_1, \dots, c_{2n-1})^T \in \mathcal{M}^0_{2n-1}$ [see Skibinsky
(1967)]. A straightforward calculation shows [see e.g.\ Szeg\"{o} (1975) or Chihara (1978)] that for a given vector
$$(c_{1}, \dots c_{2n-1})^T =
(c_1 (\eta), \dots, c_{2n-1} (\eta))^T \in \mathcal{M}^0_{2n-1}
$$
 the polynomials $P_{0,n}(x) = 1$,
\begin{equation}
\label{2.1} P_{m,n} (x) ~\colon=~ \left|\begin{array}{cccc}
  c_0& \cdots& c_{m-1}& 1\\
  \vdots& \ddots & \vdots& \vdots\\
  c_m& \cdots& c_{2m-1}& x^m
  \end{array}\right| \Big/
\left|\begin{array}{ccc}
  c_0& \cdots& c_{m-1}\\
  \vdots&  \ddots& \vdots\\
  c_ {m-1}& \cdots& c_{2m-2}
  \end{array}\right| \ ; \quad m = 1,  \dots, n
\end{equation}
have leading coefficient $1$ and are orthogonal with respect to the measure $\eta^-$, that is \be \label{2.2} \mathcal{L} (P_{k,n} \ P_{l,n}) =
\int^1_0 P_{k,n} (x) P_{l,n} (x) d \eta^- (x) = 0 \qquad \mbox{if} \quad k \neq l , \ee where the moment functional $\mathcal{L}$ is defined by
\be \label{2.3} \mathcal{L} (x^k)=  c_k  = \int^1_0 x^k d  \eta^- (x)\qquad (k = 0, \dots, 2n-1) \ . \ee
In other words: each vector $(c_1, \dots, c_{2n-1})^T \in \mathcal{M}_{2n-1}^0$ uniquely determines monic orthogonal polynomials
$P_{0,n} (x),\dots, P_{n,n}(x)$ satisfying (\ref{2.2}).
In the following we
consider a one-to-one mapping of the set $\mathcal{M}^0_{2n-1}$ on the cube $(0,1)^{2n-1}$, which was introduced by Skibinsky (1967). The new
coordinates are called canonical moments and have been studied by numerous authors [see e.g.\ Dette and Studden (1997) for a detailed
discussion]. To be precise, let $\CCp([0,1])$ denote the set of all probability measures   on the interval $[0,1]$, $ \Phi_{k}(x) =(x,\ldots
,x^{k})$ the vector of all monomials of order $k$
 and define
for a fixed vector $c =(c_1, \ldots , c_{k})^T \in {\cal M}_{k}$
$$\CCs_{k}(c):=
\left\{\mu\in\CCp([0,1]):\;\int_{0}^{1}\Phi_{k}(x)\mu(dx)=c \right\}$$ as the set of all probability measures on the interval $[0,1]$ whose
moments up to the order $k$ coincide with $c =(c_1, \ldots , c_{k})^T$. Note that this set is a singleton if and only if $c \in \partial \mathcal{M}_k$
[see Dette and Studden (1997)]. For  $k = 2,3,  \ldots$ and for a given point $(c_1, \ldots ,
c_{k-1})^T \in {\cal M}_{k-1}$ we define $c^+_k = c^+_k(c_1, \ldots , c_{k-1})$ and $c^-_k = c^-_k(c_1, \ldots , c_{k-1})$ as the largest and
smallest value of $c_k$ such that $(c_1, \ldots , c_k)^T \in \partial {\cal M}_k,$ that is
\begin{eqnarray*}
c^{-}_k &=& \min
\Bigl\{ \int^1_0 x^k\mu(dx)\mid \mu \in {\it S}_{k-1}(c_1, \ldots , c_{k-1})\Bigr\} ,\\
c^{+}_k &=& {\max} \Bigl\{ \int^1_0 x^k\mu(dx)\mid \mu \in {\it S}_{k-1}(c_1, \ldots , c_{k-1})\Bigr\}.
\end{eqnarray*}
Note that $c^-_k \le c_k \le c^+_k$ and that both inequalities are strict if and only if $(c_1, \ldots , c_{k-1})^T \in \mathcal{M}_{k-1}^0$
[see Dette and Studden (1997)].
 For a moment point $c = (c_1, \ldots , c_n)^T$ in the interior of the moment space ${\cal M}_n$ the
canonical moments or canonical coordinates of the vector $c$ are defined by
\begin{equation}\label{2.4}
p_1=c_1, \mbox{ and } p_k = \frac{c_k - c^-_k}{c^+_k - c^-_k}
 \quad k = 2, \ldots , n \ .
\end{equation}
Note that $0 < p_k < 1 $ $(k = 1, \ldots , n)$ if $(c_1,\ldots , c_n)^T \in {\cal M}_n^0$, and that the definition (\ref{2.4}) defines a one to
one mapping between $\mathcal{M}_n^0$ and the open unit cube  $(0,1)^n.$ For more details regarding canonical moments we refer to the work of
Skibinsky (1967,1968,1969) and to the monograph of Dette and Studden (1997). In particular it is shown in the lastnamed reference
 that the three term recurrence relation corresponding
to the monic orthogonal polynomials defined by (\ref{2.1}) can be represented in terms of canonical moments, that is
 $P_0 (x) = 1; P_1 (x)=x-\zeta_1~,$
{ and for } $1 \leq m \leq n-1 $
\begin{equation}\label{2.6}
P_{m+1} (x) ~=~  (x-\zeta_{2m} - \zeta_{2m+1})
  P_m (x) - \zeta_{2m-1} \zeta_{2m}
   P_{m-1} (x),
\end{equation}
where the quantities $\zeta_k$ are given by $\zeta_0=0$, $\zeta_1=p_1$, $\zeta_k =q_{k-1}p_k$ with $q_{k-1}=1-p_{k-1}$ if $k\geq 2 $. Moreover,
from the representation (\ref{2.4}) it is easy to see that $\partial c_1/ \partial p_1=1$,
\begin{eqnarray}
\frac{\partial c_k}{\partial p_j} = \left\{ \begin{array}{cll}
0 & \mbox{ if } &j > k\\
 r_k(c_1, \ldots , c_{k-1}) & \mbox{ if } &j = k\end{array} \right. \label{2.7}
\end{eqnarray}
$(k\geq 1)$ where
\begin{equation}\label{2.8}
r_{k+1} (c_1, \ldots , c_k)= r_{k+1} = c_{k+1}^+ - c_{k+1}^-= \prod^k_{j=1} p_j(1 - p_j).
\end{equation}
denotes the range of the moment space $\mathcal {M}_k$ (with the convention $r_1 = 1).$ In the following section we will use these results to
study stochastic properties of the roots of random orthogonal polynomials associated with a uniform distribution on the moment space
$\mathcal{M}_{2n-1}$.

\section{Asymptotic zero distribution of random orthogonal polynomials} 
\def\theequation{3.\arabic{equation}}
\setcounter{equation}{0}

For each $n \in \mathbb{N}$ let $C_{2n-1} = (C_{1, 2n-1}, \dots, C_{2n-1^, 2n-1})^T \sim \mathcal{U} (\mathcal{M}_{2n-1})$ denote a uniformly
distributed vector on the moment space $\mathcal{M}_{2n-1}$. From (\ref{2.0}) it follows that
$$
P \ (C_{2n-1} \in \mathcal{M}^0_{2n-1}) = 1 \ ,
$$
and consequently the random canonical moments, say $P_{1,2n-1}, \dots, P_{2n-1,2n-1}$, corresponding to the random vector $C_{2n-1}$ are well
defined with probability 1. Observing the representation (\ref{2.7}) it follows that the density of the random vector $(P_{1,2n-1}, \dots,
P_{2n-1,2n-1})$ is given by
\begin{equation}\label{3.1}
\prod^{2n-2}_{j=1} \frac{\Gamma(4n-2j)}{(\Gamma(2n-j))^2} \prod^{2n-2}_{j=1} (p_i (1 - p_i))^{2n-j-1},
\end{equation}
 which means that $\{
(P_{j,2n-1})^{2n-1}_{j=1} \ | \ n \in \mathbb{N} \}$ is a triangular array of rowwise independent random variables, where $P_{j,2n-1}$ has a
symmetric Beta-distribution on the interval $[0,1]$ with parameter $2n-j$, that is $P_{j,2n-1} \sim B (2n-j, 2n-j)$. In what follows, let
$P_{m,n} (x)$ denote the $m$th random monic orthogonal polynomial associated with the random vector $C_{2n-1}$ by equation (\ref{2.1}) and $X_{1,n},
\dots, X_{m,n}$ the corresponding roots which are real with probability 1 [see Szeg\"{o} (1975)]. Our first result gives an explicit representation
for the joint density of the
random vector $(X_{1,n}, \dots, X_{n,n})^T$.

\bigskip

{\bf Theorem 3.1.} {\it The joint density of the roots $X_{1,n}, \dots, X_{n,n}$ of the monic random orthogonal polynomial $P_{n,n} (x)$
corresponding to a random vector $C_{2n-1} = (C_{1,2n-1}, \dots, C_{2n-1,2n-1})^T \sim \mathcal{U} (\mathcal{M}_{2n-1})$ by equation
(\ref{2.1}) is given by \be \label{3.2a} f (x_1, \dots, x_n) = c \cdot \prod_{1 \leq i < j \leq n} | x_i - x_j |^4 \prod^n_{i=1} I_{[0,1]}
(x_i), \ee where the normalizing constant $c$ is defined by
$$
c = {1\over \Gamma (n+1) } \prod^{n-1}_{r=0} \frac {\Gamma (2r + 2n)}{\Gamma (2r + 1)^2 \Gamma (2r+2)}.
$$
}

\bigskip

{\bf Proof.} Consider the random canonical moments $(P_{1,2n-1}, \dots, P_{2n-1,2n-1})^T$ corresponding to the random vector $(C_{1,2n-1,} \dots,
C_{2n-1,2n-1})^T \sim \mathcal{U} (\mathcal{M}_{2n-1})$. Let $\Xi_1 = \Xi_{1,2n-1} = P_{1,2n-1}$ and for $j\geq 2$ $\Xi_j =\Xi_{j,2n-1} =
(1-P_{j-1,2n-1}) P_{j,2n-1}$, then it follows from the recursive relation (\ref{2.6}) that the random orthogonal polynomial $P_{n,n} (x)$
 can be represented as the determinant of a symmetric tridiagonal matrix, that is
 \begin{eqnarray}
 \label{3.2}
~~~~~~~~~~~~ P_{n,n} (x) = \left | \begin{array}{cccccc}
{ x - \Xi_{1}} &  - { \sqrt{\Xi_{1} \Xi_{2}}} & & & & \\
- { \sqrt{\Xi_{1} \Xi_{2}}} & { x - \Xi_{2} - \Xi_{3} } & & & & \\
& - { \sqrt{\Xi_{3} \Xi_{4}} } & \ddots &  & & \\
& & & &\ddots & \\
& & \ddots & & & - { \sqrt{\Xi_{2n-3} \Xi_{2n-2}}} \\
& & & & - { \sqrt{\Xi_{2n-3} \Xi_{2n-2}}} & { x -  \Xi_{2n-2}- \Xi_{2n-1}}
\end{array} \right |
\end{eqnarray}
Consequently, the roots $X_{1,n}, \dots, X_{n,n}$ of the polynomial $P_{n,n} (x)$ are the eigenvalues of the random Jacobi matrix $J_{n,n}$
where for $m<n$
\begin{eqnarray}
\label{3.3} ~~~~~~~~~~~~ J_{m,n} = \left ( \begin{array}{cccccc}
{  \Xi_{1}} &   { \sqrt{\Xi_{1} \Xi_{2}}} & & &  &\\
 { \sqrt{\Xi_{1} \Xi_{2}}} & {  \Xi_{2} + \Xi_{3} } &  & & & \\
&  { \sqrt{\Xi_{3} \Xi_{4}} } & \ddots & & & \\
& & & & \ddots& \\
& & \ddots & & &  { \sqrt{\Xi_{2m-3} \Xi_{2m-2}}}\\
& & & &  { \sqrt{\Xi_{2m-3} \Xi_{2m-2}}} & {  \Xi_{2m-2}+ \Xi_{2m-1}}
\end{array} \right )
\in \mathbb{R}^{m \times m} \ .
\end{eqnarray}
Now define for $j = 0, \dots, 2n-1$ the quantities $\alpha_j = 2P_{j+1,2n-1} - 1$ then $\alpha_0, \dots, \alpha_{2n-1}$ are independent random
variables and $\alpha_j$ follows a symmetric Beta-distribution on the interval $[-1,1]$ with parameter $2n-j-1$. Moreover, a straightforward
calculation shows that the matrix $\tilde J_n = 4J_{n,n}-I_n$ can be represented as \be \label{3.4} \tilde J_n := \begin{pmatrix}
                b_1 & a_1    &         &         \\
                a_1 & b_2    & \ddots  &         \\
                    & \ddots & \ddots  & a_{n-1} \\
                    &        & a_{n-1} & b_n
        \end{pmatrix} \in \mathbb{R}^{n \times n}
 \ee
with entries
\begin{eqnarray*}
b_{k+1} &=& (1-\alpha _{2k-1} )\alpha _{2k} - (1+\alpha _{2k-1} )\alpha _{2k-2} \\
a_{k+1} &=& \left \{ (1- \alpha _{2k-1} )(1-\alpha _{2k} ^2)(1+\alpha _{2k+1} ) \right \} ^{1/2} \ .
\end{eqnarray*}
and $\alpha_{-1}=-1$. This matrix has been considered recently by Killip and Nenciu (2004) in a more general context. By the results of these
authors we obtain that the density of the random eigenvalues $\lambda_1, \dots, \lambda_n$ of the matrix $\tilde J_n$ is given by
$\tilde c \prod_{1 \leq i < j \leq n} | \lambda_i - \lambda_j |^4$, where $\tilde c$ is an appropriate normalizing constant. 
Transferring this result to the eigenvalues of the matrix $J_{n,n}$ it follows that the density of the roots of the random polynomial $P_{n,n}
(x)$ is given by (\ref{3.2a}). \hfill $\Box$

\bigskip \bigskip

The density defined by (\ref{3.2a}) is a special case of the Jacobi ensemble in the symplectic case [see e.g.\ Mehta (2004)] which has found considerable
interest in the recent literature  [see e.g.\ Collins (2005) or Johnstone (2008) among others]. The density of the general Jacobi
$\beta$-ensemble is given by
$$
c^n_{a_n,b_n}  \prod_{1 \leq i \leq j \leq n} | \lambda_i - \lambda_j |^\beta \: \prod^n_{i=1} \lambda_i^{a_n} (1-
\lambda_i)^{b_n} \: I_{[0,1]} (\lambda_i)
$$
where $a_n,b_n \ge -1$, $c^n_{a_n, b_n}$ is a normalizing constant and $\beta > 0$ [see Killip and Nenciu (2004)]. While most authors consider the case where $n
\to \infty$, $a_n/n \to \alpha$; $b_n/n \to \beta$, much less attention has been paid to the case where $a_n$ and $b_n$ are fixed.  We use the
explicit representation for the density of the roots of $P_{n,n} (x)$ in terms of eigenvalues of the random matrix $J_{n,n}$ defined in
(\ref{3.3}) to derive asymptotic properties for the empirical distribution function of the random variables $X_{1,n},\dots, X_{n,n}$.

\bigskip

{\bf Theorem 3.2.} {\it Let
$$
F_n (x) = \frac {1}{n} \sum^n_{j=1} I \{ X_{j,n} \leq x \}
$$
denote the empirical distribution function of the roots of the monic random orthogonal polynomial $P_{n,n} (x)$ associated with the random
vector $C_{2n-1} = (C_{1,2n-1}, \dots, C_{2n-1,2n-1})^T \sim \mathcal{U}(\mathcal{M}_{2n-1})$, then for all $x \in [0,1]$
$$
\lim_{n \to \infty} F_n (x) = \frac {1}{\pi} \int^x_0 \frac {dt}{\sqrt{t (1-t)}} \qquad \mbox{a.s.}
$$}

\bigskip

{\bf Proof.} Define
\begin{eqnarray}\label{3.5}
D_{m,n} = \left (
\begin{array}{lllll}
\frac {1}{2} & \frac {1}{2\sqrt{2}} &&& \\
\frac {1}{2\sqrt{2}} & \frac {1}{2} & \frac {1}{4} & & \\
& \frac {1}{4} & \frac {1}{2} & & \\
& & & \ddots & \frac {1}{4} \\
& & & \frac {1}{4} & \frac {1}{2}
\end{array}
\right ) \in \mathbb{R}^{m \times m},
\end{eqnarray}
then a straightforward calculation shows that the characteristic polynomial of the matrix $D_{n,n}$ is given by the monic Chebyshev polynomial
of the first kind on the interval $[0,1]$, that is \be \label{ch1} \overline T_n (x) : = \det (x I_n - D_{n,n}) = 2^{-2n+1} \cos (n \arccos
(2x-1)) , \ee which has roots \be \label{3.6} x_{k,n} = \frac {\cos ( (2k-1) \frac {\pi}{2n}) + 1}{2} \ ; \qquad k=1,\dots, n . \ee Moreover,
define $T_k (x) = 2^{2k-1} \ \overline T_k (x)$ for  $k\geq1$, $T_0 (x) = 1/\sqrt{2}$ and $T (x) = (T_0 (x), \dots, T_{n-1} (x) )^T$, then it
follows from the relation $\cos ( (n+1) z) + \cos ( (n-1) z) = 2 \cos (nz) \cos (z)$ by a straightforward calculation that \be \label{eig}
x_{k,n} \ T (x_{k,n}) = D_{n,n} \ T (x_{k,n}) \ ; \qquad k = 1, \dots, n \ . \ee This shows that the vectors
\begin{eqnarray} \label{eigvec} ~~~~~~~~
t_{k,n} = T
(x_{k,n}) = \left ( 1 / \sqrt{2},  \cos ( \frac {2k-1}{2n} \pi), \dots, \cos ((n-1) \frac {2k-1}{2n} \pi) \right )^T; ~ k = 1, \dots, n
\end{eqnarray}
are the eigenvectors of the matrix $D_{n,n}$ corresponding to the eigenvalues $x_{1,n}, \dots, x_{n,n}$, respectively. In the following
discussion define
$$
N_n (x) = \frac {1}{n} \sum^n_{j=1} I \{ x_{j,n} \leq x \}
$$
as the empirical distribution function of the roots of $\overline T_n (x)$ (or $T_n (x)$), then an elementary calculation shows  that for all
$x \in [0,1]$ \be \label{3.7} \lim_{n \to \infty} N_n (x) = \frac {1}{\pi} \int^x_0 \frac {dt}{\sqrt{t (1-t)}}. \ee From Bai (1999) we have for
the Levy distance between the empirical distribution functions $F_n$ and $N_n$ the estimate
\begin{eqnarray} \label{3.8} \nonumber
L^3 (F_n, N_n) & \leq & \frac {1}{n} \sum^n_{j=1} \mid X_{(j),n} - x_{(j),n} \mid^2 \ \leq \frac {1}{n} \ tr \ (J_{n,n}-D_{n,n})^2 \\
& \leq & \frac {1}{n}  \sum^n_{i=1} (J_{ii} - D_{ii})^2 + \frac {2}{n} \sum^{n-1}_{i=1} (J_{ii+1} - D_{ii+1})^2,
\end{eqnarray}
where $X_{(1),n} \leq \dots \leq X_{(n),n}$ and  $x_{(1),n} \leq \dots \leq x_{(n),n}$ denote the ordered roots of the polynomials $P_{n,n}
(x)$ and $\overline T_n (x)$, respectively, and $J_{ij}$ and $D_{ij}$ are the elements of the tridiagonal matrices $J_{n,n}$ and $D_{n,n}$
defined in (\ref{3.3}) and (\ref{3.5}), respectively. Using the notation $Q_{j,2n-1} = 1 - P_{j,2n-1}$ $(j=1, \dots, 2n-1)$ it now follows by a
 straightforward calculation that for $2 \leq i \leq n$
\begin{eqnarray*}
(J_{ii} - D_{ii})^2 &=& \biggl(Q_{2i-3, 2n-1} P_{2i-2,2n-1} + Q_{2i-2, 2n-1} P_{2i-1, 2n-1} - \frac {1}{2}\biggr)^2 \\
& \leq & \biggl \{ \big | Q_{2i-3,2n-1} - \frac {1}{2} \big | + \frac {1}{2} \big | P_{2i-2,2n-1} - \frac {1}{2} \big | + \big | Q_{2i-2,2n-1}
-
\frac {1}{2} \big | + \frac {1}{2} \big | P_{2i-1,2n-1} - \frac {1}{2} \big |  \biggr \}^2 \\
& \leq & \big | P_{2i-3,2n-1} - \frac {1}{2} \big |^2 + 9 \big | P_{2i-2,2n-1} - \frac {1}{2} \big |^2 + \big | P_{2i-1,2n-1} -\frac {1}{2}
\big |^2
\end{eqnarray*}
and $(J_{11} - D_{11})^2 \leq | P_{1,2n-1} - \frac {1}{2} |^2$. Similarly we obtain for $2  \leq i \leq n-1$
\begin{eqnarray*}
\big | J_{ii+1} - D_{ii+1}\big |^2 &\leq& \big |Q_{2i-2, 2n-1} P_{2i-1,2n-1}  Q_{2i-1, 2n-1} P_{2i, 2n-1} - \frac {1}{16}\big | \\
& \leq & \frac {1}{8} \biggl \{  2 \big | P_{2i-2,2n-1} - \frac {1}{2} \big | + 3  \big | P_{2i-1,2n-1} - \frac {1}{2} \big | + \big |
P_{2i,2n-1} - \frac {1}{2} \big |   \biggr \}
\end{eqnarray*}
and $\big | J_{1,2} - D_{1,2}\big |^2 \leq  \frac {5}{2} \big |P_{1,2n-1} - \frac {1}{2} \big | + \frac {1}{4} \big |P_{2,2n-1} - \frac {1}{2}
\big | $. In the following discussion we will show that
\begin{eqnarray} \label{3.10}
&& \frac {1}{n} \sum^{n-1}_{i=1} \  \big | P_{2i,2n-1} - \frac {1}{2}  \big |^2 \ \ \    \stackrel {\scriptstyle a.s.}  {\longrightarrow} \ 0 ,
\\ \label{3.11}
&& \frac {1}{n} \sum^{n}_{i=1} \  \big | P_{2i-1,2n-1} - \frac {1}{2}  \big |^2 \ \stackrel {\scriptstyle a.s.}  {\longrightarrow} \ 0 ,
\end{eqnarray}
which directly implies
$$
\frac {1}{n} \sum^{n}_{i=1} \  (J_{ii} - D_{ii})^2 \ \stackrel {\scriptstyle a.s.}  {\longrightarrow} \ 0 .
$$
Moreover, for the remaining sum in (\ref{3.8}) it follows that
\begin{eqnarray*}
\frac {1}{n} \sum^{n-1}_{i=1} (J_{ii+1} - D_{ii+1})^2 & \leq & \gamma \biggl \{ \frac {1}{n} \sum^{n-1}_{i=1} \big | P_{2i,2n-1} - \frac {1}{2} \big
|^2 \biggr \}^{1/2} + \biggl \{ \frac {1}{n} \sum^{n}_{i=1} \big | P_{2i-1,2n-1} -
\frac {1}{2} \big |^2 \biggr \}^{1/2} \\
& \stackrel {\scriptstyle a.s.} {=} & o \ (1),
\end{eqnarray*}
where the constant $\gamma$ does not depend on $n$. The assertion is now a consequence of (\ref{3.7}) and (\ref{3.8}), which yields for the Levy-distance between the empirical distribution
function $F_n$ and distribution function $F$ of the arcsine measure
$$
L (F_n, F) \leq L (F_n, N_n) + L (N_n,F) = o \ (1) \quad \mbox{a.s.}
$$

For a proof of the almost sure convergence in (\ref{3.10}) and (\ref{3.11}) we restrict ourselves to the statement (\ref{3.10}), the remaining
case is treated similarly. In order to prove  (\ref{3.10}) we will use a strong law of large numbers for arrays of rowwise independent random
variables. To be precise, define
$$
Z_{n,k} = (P_{2k,2n-1} - \frac {1}{2})^2 - \frac {1}{8 (2n-2k)+4} \qquad k=1,\dots,n-1
$$
then $E[Z_{n,k}]=0$ and a  tedious calculation shows that for $k = 1, \dots, n-1$ and any $t > 0$
$$
P ( | Z_{n,k} | > t) \leq P ( | Z_{n,n-1} | > t) \ ,
$$
where the distribution of the random variable $Z_{n,n-1}$ does not depend on $n$ [note that $P_{2n-2,2n-1} \sim B(2,2)$]. Consequently, we
obtain from Theorem 2 in Hu, M\'{o}ricz and Taylor (1989) that
$$
Z_n = \frac {1}{n^{1/p}} \sum^{n-1}_{k=1} Z_{n,k} \  \stackrel {\scriptstyle{a.s.}}{\longrightarrow} \ 0
$$
for any $p \in [1,2)$. Observing that
$$
\frac {1}{4} \sum^{n-1}_{k=1} \frac {1}{2(2n-1-2k)+1} = O (\log n)
$$
it follows for any $p \in [1,2)$ that
$$
\frac {1}{n} \sum^{n-1}_{k=1} \big | P_{2k,2n-1} - \frac {1}{2} \big |^2 = \frac {1}{n^{1-1/p}} \  \Big \{ Z_n + O \ \Big(\frac {\log
n}{n^{1/p}}\Big) \Big \} = o \left ( \frac {1}{n^{1 - 1/p} }\right ) \: \: \mbox{a.s.},
$$
which establishes (\ref{3.10}) and completes the proof of Theorem 3.2. \hfill $\Box$

\bigskip





%
%
%
%

\bigskip

Our final result refers to the asymptotic behaviour of the roots of the $m$th orthogonal polynomial associated with the random moment vector
$C_{2n-1} = (C_{1,2n-1,} \dots, C_{2n-1,2n-1})^T \sim \mathcal{U} (\mathcal{M}_{2n-1})$, where $m$ is fixed. In this case, the limiting
distribution is normal, where the roots of the Chebyshev polynomial $\overline T_n (x)$ defined in (\ref{ch1}) are used for the centering.

\bigskip

{\bf Theorem 3.3.} {\it Let $X_{1,n}, \dots, X_{m,n}$ denote the roots of the $m$-th random monic orthogonal polynomial $P_{m,n} (x)$, which corresponds to the
random vector $C_{2n-1} \sim \mathcal{U} (\mathcal{M}_{2n-1})$ by equation (\ref{2.1}), then
$$
4 \sqrt{n} \left \{ (X_{1,n}, \dots, X_{m,n})^T - (x_{1,n}, \dots, x_{m,n})^T \right \} \textstyle{ {\mathcal{D} \atop \overrightarrow{n \to
\infty}}} \ \mathcal{N} (0, \Gamma),
$$
where $x_{1,m}, \dots, x_{m,m}$ are the roots of the Chebyshev polynomial of the first kind defined by (\ref{3.6}) and the matrix $\Gamma$ is
given by $\Gamma = \frac {2}{m} (\gamma_{k,l})^m_{k,l=1}$, where
\begin{eqnarray*}
\gamma_{k,l} &=& \frac {1}{4} + \frac {1}{2} \sum^m_{j=2} \ T^2_{j-1} (x_{k,m}) \ T^2_{j-1} (x_{l,m}) - \frac {1}{4} \sum^{m-1}_{j=1} \
\left ( T^2_{j-1} (x_{k,m}) \ T^2_j (x_{l,m}) + T^2_{j-1} (x_{l,m}) \ T^2_{j} (x_{k,m} ) \right ) \\
& + & \frac {1}{4} \  T_1  (x_{k,m}) \ T_1 (x_{l,m}) + \frac {1}{2} \sum^{m-1}_{j=2} \
T_{j-1} (x_{k,m}) \ T_j (x_{k,m}) \ T_{j-1} (x_{l,m}) \ T_{j} (x_{l,m}) \\
& - & \frac {1}{4} \sum^{m-2}_{i=1} \left ( T_{j-1} (x_{k,m}) \ T_j (x_{k,m}) \ T_j (x_{l,m}) \ T_{j+1} (x_{l,m}) + T_{j-1} (x_{l,m}) \ T_j
(x_{l,m}) \ T_j (x_{k,m}) \ T_{j+1} (x_{k,m}) \right ),
\end{eqnarray*}
 $T_j (x)= \cos ( j \arccos (2x-1))$ and $x_{k,m}$ is the $k$th zero of the polynomial $T_m (x)$ defined in (\ref{3.6}).}

\bigskip

{\bf Proof.} It is easy to see that for fixed $m \in \mathbb{N}$ we have for $k = 1, \dots, m$ \be \label{3.13} 4 \sqrt{n} \ \Big(P_{k,2n-1} -
\frac {1}{2}\Big)  \ \textstyle{ {\mathcal{D} \atop \overrightarrow{n \to \infty}}} \ \mathcal{N} (0,1) \ , \ee which implies that the vector
$4 \sqrt{n} (P_{1,2n-1}, \dots,  P_{m,2n-1})^T$ is asymptotically multivariate normal distributed with mean 0 and covariance matrix $I_m$ (note
that $P_{1,2n-1},\dots,P_{m,2n-1}$ are independent random variables). A straightforward calculation now shows that \be \label{3.14} 4 \sqrt{n} \
\Big \{  \Big (\Xi_{1,2n-1}, \dots, \Xi_{2m-1,2n-1} \Big)^T - \Big(\frac {1}{2}, \frac {1}{4}, \dots, \frac {1}{4}\Big)^T \Big  \} \textstyle{
{\mathcal{D} \atop \overrightarrow{n \to \infty}}} \: \mathcal{N}_{2m-1} (0,A_{2m-1}) \ , \ee where the matrix $A_{2m-1} \in \mathbb{R}^{2m-1
\times 2m-1}$ is tridiagonal and given by
\begin{eqnarray} \label{3.15}
A_{2m-1} = \left (
\begin{array}{rrrrrr}
1 &  - \frac {1}{2} &&&& \\
- \frac {1}{2} & \frac {1}{2} & - \frac {1}{4} &&& \\
& -\frac {1}{4} & \frac {1}{2} & - \frac {1}{4} & & \\
&&\ddots&\ddots & \ddots & \\
&&&&\ddots &- \frac {1}{4} \\
&&&& - \frac {1}{4} & \frac {1}{2}
\end{array}
\right ).
\end{eqnarray}
A further application of the delta-method finally yields
\begin{eqnarray} \label{3.16}
4 \sqrt{n} \  \Bigg \{ \left (
\begin{array}{c}
\Xi_{1} \\ \Xi_{2} + \Xi_{3} \\  \vdots \\
\Xi_{2m-2} + \Xi_{2m-1} \\
\sqrt{\Xi_{1}\Xi_{2}} \\
\sqrt{\Xi_{3}\Xi_{4}} \\
\vdots \\
\sqrt{\Xi_{2m-3}\Xi_{2m-2}}
\end{array}
 \right ) - \left (
\begin{array}{c}
 {1}/{2} \\
 {1}/{2} \\ \vdots \\  {1}/{2} \\  {1}/({2} \sqrt{2}) \\  {1}/{2} \\ \vdots \\{1}/{2}
\end{array}
 \right ) \Bigg \}
 \textstyle{ {\mathcal{D} \atop \overrightarrow{n \to \infty}}} \ \mathcal{N} (0,V)
\end{eqnarray}
where the covariance matrix $V = diag (V_{11}, V_{22}) \in \mathbb{R}^{2m-1 \times 2m-1}$ is block diagonal with blocks $V_{11} = A_m \in
\mathbb{R}^{m \times m}$ and
\begin{eqnarray*}
V_{22} = \frac {1}{8} \left (
\begin{array}{rrrrrr}
1 &  - \frac {1}{\sqrt{2}} &&&& \\
- \frac {1}{\sqrt{2}} & 1 & - \frac {1}{2} &&& \\
& -\frac {1}{2} & 1 & - \frac {1}{2} & & \\
&&\ddots&\ddots & \ddots & \\
&&&&\ddots &- \frac {1}{2} \\
&&&& - \frac {1}{2} & 1
\end{array}
\right ) \in \mathbb{R}^{m-1 \times m-1}  .
\end{eqnarray*}
Consequently, the difference of the matrices $J_{m,n} - D_{m,n}$ defined by (\ref{3.3}) and (\ref{3.5}) converges also weakly, that is
$$
\lim_{n \to \infty} 4 \sqrt{n} \ (J_{m,n} - D_m) = \mathcal{S} \qquad \mbox{(weakly)},
$$
where $\mathcal{S}$ denotes a triangular matrix given by
\begin{eqnarray*}
\mathcal{S}  = \left (
\begin{array}{lllll}
M_1 & N_1 &&& \\
N_1 & M_2 & N_2 & & \\
& N_2 & M_3 & N_3 & \\
&& \ddots & \ddots & \ddots \\
&&& \ddots & N_{m-1} \\
&&& N_{m-1} & M_{m}
\end{array}
\right ) \ ,
\end{eqnarray*}
and $\{ M_1, \dots, M_m  \}$ and $\{ N_1, \dots, N_{m-1} \}$ define independent samples $M_1 \sim \mathcal{N}(0,1)$, $M_j \sim \mathcal{N}
(0,1/2)$ $(j = 2, \dots, m)$ $N_j \sim \mathcal{N} (0,1/8)$ $(j = 1, \dots, m-1)$ such that
\begin{eqnarray}\label{c1}
\rm{Cov} \ (M_i, M_j) &=& \left \{
\begin{array}{lllrl}
-\frac {1}{2} &   \mbox{if} &    i=1,j=2  & \mbox{or} &   i=2,j=1 \\
-\frac {1}{4} &    \mbox{if}&    j=i+1     & \mbox{or} &   j=i-1 \\
\: \: \: \: 0 & \mbox{else}
\end{array}
\right .\\  \label{c2} \rm{Cov} \ (N_i, N_j) &=& \left \{
\begin{array}{lllrl}
-\frac {1}{8\sqrt{2}} &   \mbox{if} &    i=1,j=2  & \mbox{or} &   i=2,j=1 \\
-\frac {1}{16} &    \mbox{if}&    j=i+1     & \mbox{or} &   j=i-1 \\
\: \: \: \: 0& \mbox{else}
\end{array}
\right .
\end{eqnarray}

Using the same arguments as in Dimitriu and Edelman (2005) it now follows that the first order properties of the matrix $J_{m,n}$ are the same
as the first order properties of the matrix $D_{m,n} + \frac {1}{4 \sqrt{n}}\mathcal S$. In particular, using Lemma 2.1 in this reference we
obtain for $i = 1, \dots, m$ that the asymptotic distribution of the vector
$$
4 \sqrt{n} \left \{  (X_{1,n}, \dots, X_{m,n})^T - (x_{1,n}, \dots x_{m,n} )^T \right \}
$$
coincides with the  distribution of the random vector \be \label{3.17} G = \left (  \frac {t^T_{1,m} \ \mathcal{S} \ t_{1,m}}{t^T_{1,m}
t_{1,m}}, \dots, \frac {t^T_{m,m} \ \mathcal{S} \ t_{m,m}}{t^T_{m,m} t_{m,m}}  \right )^T \ee where the vectors $t_{1,m}, \dots, t_{m,m} \in
\mathbb{R}^m$ are defined by (\ref{eigvec}). Obviously, $G$ has a multivariate normal distribution with mean zero and it remains to
calculate the corresponding covariance matrix. For this purpose we note that for $k=1,\dots,m$
$$
t^T_{k,m} \ t_{k,m} = \frac {1}{2} + \sum^{m-1}_{l=1} \cos^2 (l \ \frac {2k-1}{2m} \pi) = \frac{m}{2} \ .
$$
Therefore it remains to calculate the covariance matrix of the vector $(t^T_{1,m} \ \mathcal S \ t_{1,m}, \dots,  t^T_{m,m} \ \mathcal S \
t_{m,m})$. For this purpose we denote by $(v_1, \dots, v_m)$ and $(w_1, \dots, w_m)$ the components of the vectors $t_{k,m}$ and $t_{l,m}$,
respectively, and obtain
\begin{eqnarray*}
&& {\rm Cov} \ (t^T_{k,m} \ \mathcal S \ t_{k,m}, \ t^T_{l,m} \ \mathcal S \ t_{l,m}) = {\rm Cov} \Big ( \sum^m_{i=1} v^2_i M_i, \sum^m_{j=1}
w^2_j M_j \Big )
+ 4 \ {\rm Cov} \Big ( \sum^{m-1}_{i=1}  v_i v_{i+1} N_i, \sum^{m-1}_{j=1} w_j w_{j+1} N_j \Big ) \\
&&= v^2_1 w^2_1 + \frac {1}{2} \sum^m_{j=2} v^2_j w^2_j - \frac {1}{4} \Big \{ w^2_2 + v^2_2 +
\sum^{m-1}_{i=2} (v^2_i w^2_{i+1} + v^2_{i+1} w^2_i )  \Big \} \\
&&+ \frac {1}{2} \sum^{m-1}_{i=1} v_i v_{i+1} w_i w_{i+1} - \frac {1}{4} \Big \{v_2 w_2 w_3 + w_2 v_2 v_3 + \sum^{m-2}_{i=2} (v_i v_{i+1}
w_{i+1} w_{i+2} + v_{i+1} v_{i+2} w_i w_{i+1} ) \Big \},
\end{eqnarray*}
where we have used (\ref{c1}) and (\ref{c2}) and $w_1 = v_1 = 1/\sqrt{2}$.
This proves the assertion of the Theorem. \hfill $\Box$.

\medskip
\bigskip

{\bf Acknowledgements.} The authors are grateful to Martina Stein  who typed most parts of this paper with considerable technical expertise.
The work of the authors was supported by the Sonderforschungsbereich Tr/12, Fluctuations and universality of invariant random matrix ensembles
(project C2) and in part by a NIH grant award IR01GM072876:01A1.

\bigskip
\bigskip

{\Large References}

\bigskip

Z. D. Bai (1999).  Methodologies in spectral analysis of large dimensional random matrices, a review. Statistica Sinica 9, 611-677.

\smallskip

{F. C. Chang, J. H. B. Kemperman and W. J. Studden (1993). A normal limit theorem for moment sequences. Ann. Probab. 21, 1295-1309.}

\smallskip

{T. Chihara (1978). An Introduction to Orthogonal Polynomials. Gordon and Breach, New York.}

\smallskip

H. Dette, F. Gamboa (2007). Asymptotic properties of the algebraic moment range process. Acta Math.\ Hung.\ 116, 247-264.

\smallskip

H. Dette, W.~J.  Studden (1997).
 { The theory of canonical moments with applications in statistics,
probability, and analysis}. John Wiley \& Sons Inc., New York.
 ISBN 0-471-10991-6.
 A Wiley-Interscience Publication.

\smallskip

B. Collins (2005). Product of random projections, Jacobi ensembles and universality problems arising from free probability.
    Probability Theory and Related Fields 133, 315-344.

 \smallskip

 J. Dumitriu, A. Edelman (2005). Eigenvalues of Hermite and Laguerre ensembles: large beta asymptotics.
 Ann.\ I.\ Henri Poincare, Probabilit\'{e}s et Statistiques 41, 1083-1099.

\smallskip

F. Gamboa and L.V. Lozada-Chang (2004). Large deviations for random power moment problem. Ann. Probab. 32, 2819-2837.

\smallskip

T.-C. Hu, F. M\'{o}ricz, R. L. Taylor (1989). Strong laws of large numbers for arrays of rowwise independent random variables. Acta Math.\
Hung. 54, 153-162.

\smallskip

I.~M. Johnstone (2008). Multivariate analysis and Jacobi ensembles: largest eigenvalue, Tracy Widom limits and rates of convergence. Annals of
Statistics, to appear.

\smallskip

S. Karlin, L.~S.~Shapeley (1953).
 Geometry of moment spaces. Amer. Math. Soc. Memoir No. 12, Amer. Math. Soc., Providence, Rhode Island.

\smallskip

{R. Killip, I. Nenciu I (2004). Matrix models for circular ensembles. Int. Math. Res. Not. 50, 2665-701.}
\smallskip


\smallskip

M. L. Mehta (2004). Random Matrices. Academic Press.
\smallskip

M. Skibinsky (1967).
 The range of the $(n+1)$-th moment for distributions on $[0,\,1]$.
{ J. Appl. Probability} 4, 543--552.

\smallskip

M. Skibinsky (1968).
 Extreme $n$th moments for distributions on $[0,\,1]$ and the inverse
of a moment space map. J. Appl. Probability 5,  693--701.

\smallskip
M. Skibinsky (1969).
 Some striking properties of binomial and beta moments. {Ann. Math. Statist.} 40,  1753--1764.

\smallskip

{G. Szeg\"{o} (1975). Orthogonal Polynomials. Amer. Math. Soc. Colloqu. Publ. 23, Providence, RI. }

\end{document}